\documentclass[12pt]{article}    
\usepackage{amsmath}
\usepackage{amssymb}
\usepackage{amsmath,amssymb,amsfonts,amsthm,latexsym}
\usepackage{booktabs}
\usepackage{array} 
\usepackage{lscape}
\usepackage{graphicx}
\usepackage[all]{xy}
\usepackage{epsfig}
\usepackage{colordvi}
\usepackage{latexsym}
\usepackage{authblk}
\def\mysection{\setcounter{equation}{0}\section}
\newtheorem{Lemma}{Lemma}[section]
\newtheorem{Theorem}[Lemma]{Theorem}
\newtheorem{Proposition}[Lemma]{Proposition}

\theoremstyle{definition}
\newtheorem{defn}{Definition}[section]

\newcommand{\beq}{\begin{equation}}
\newcommand{\eeq}{\end{equation}}
\newcommand{\beqr}{\begin{eqnarray}}
\newcommand{\eeqr}{\end{eqnarray}}

\begin{document}
	\title{\Large {{Topogenous structures on faithful and amnestic functors} }}
    \author[$*,\;a$]{Minani Iragi}
    \author[$b$]{David Holgate}
	\author[$c$]{Josef Slapal}
	
	\affil[$^{a,\;c}$]{{\footnotesize Institute of Mathematics, Brno University of Technology, Technick\'{a} 2, 616 69, Brno, Czech Republic.}}
	
	\affil[$b$]{{\footnotesize Department of Mathematics and Applied Mathematics\\ University of the Western Cape,
			Bellville 7535, South Africa.}}
	
	\date{}
	\maketitle
	\def\thefootnote{\fnsymbol{footnote}}
	\setcounter{footnote}{0}
	
	\footnotetext{ 
	E-mail addresses: $^{a}$Minani.Iragi@vutbr.cz, $^{b}$dholgate@uwc.ac.za, $^{c}$slapal@fme.vutbr.cz.\\
	$^{*}$Corresponding author.\\	
	The first author acknowledges the support from the Brno University of Technology (BUT) under the project MeMoV II no. CZ.02.2.69/0.0/0.0/18-053/0016962. The second author acknowledges the National Research Foundation of South Africa. The third author acknowledges the support by BUT from the Specific Research Project no. FSI-S-20-6187.}
	\begin{abstract}
Departing from a suitable categorical concept of topogenous orders defined relative to the bifibration of subobjects, this note introduces and studies topogenous orders on faithful and amnestic functors. Amongst other things, it is shown that this approach captures the formal closure operators and leads to the introduction of formal interior operators. Turning to special morphisms relative to the orders introduced, we show that a morphism is strict relative to an order if the order preserves codomains of its cocartesian liftings  while a morphism is final if the order reflects domains of its cartesian liftings.  Key examples in topology and algebra that demonstrate our results are included. 	
    \end{abstract}
\\

\mysection{Introduction}
{\indent A topogenous order $\sqsubset$ on a category $\mathcal{C}$ equipped with a proper ($\mathcal{E}$, $\mathcal{M}$)-factorization structure for morphisms is a family of binary relations, each on the subobject lattice, sub$X$, for an object $X$ in $\mathcal{C}$ (subject to some axioms) (\cite{MR0286814}). This notion, which is crutial for the syntopogenous structures introduced by Cs{\'a}sz{\'a}r (\cite{MR0157340}) with the aim of proposing a unified approach to topological, uniform and proximity spaces, has played a salient role in providing a single setting study of categorical closure (\cite{dikranjan1987closure}), interior (\cite{MR2838385}) and neighbourhood (\cite{MR1810290}) operators and led to the introduction of quasi-uniform structures in categories (see e.g \cite{holgate2019quasi, MR4239728, iragi2019quasi, MR3641225, MR4741225}). Topogenous orders are easier to work with when it comes to the study of topological structures on categories.}

{\indent  A new way of introducing topological structures on categories was proposed by Dikranjan and Giuli who originally defined the categorical closure operators. A closure operator on a finitely $\mathcal{M}$-complete category $\mathcal{C}$ is a pointed endofunctor of $\mathcal{M}$, where the class $\mathcal{M}$ is seen as the full subcategory of the arrow category $\mathcal{C}^{2}$ whose objects are the morphisms from $\mathcal{M}$, which “commutes” with  the codomain functor  $ cod: \mathcal{M}\longrightarrow \mathcal{C}$. The functor $cod$ is a bifibration and indeed, categorical closure operators are defined relative to the bifibration of $\mathcal{M}$-subobjects. Replacing $\mathcal{M}$ by an arbitrary category and the codomain functor by a faithul and amnestic functor $F$, called $form$ in \cite{MR3226619}, the notion of formal closure  operator was recently obtained in \cite{MR3641248}. The formal closure operators behave well and keep most of the properties of the categorical closure operators. In particular, these closure operators capture epireflective subcategories through the notions of idempotency, coheredity, and minimality. 

 {\indent  Viewing a categorical topogenous order as defined relative to the bifibration of subobjects, we replace the functor $ cod: \mathcal{M}\longrightarrow \mathcal{C}$ by a suitable form over $\mathcal{C}$, described at the begining of the third section of this note. In this case, the subobject lattice of an object $X\in \mathcal{C}$ is thought of as the $fibre$ over $X$ so that for any $f: X\longrightarrow Y$ in $\mathcal{C}$, we think of the image of a subobject as the codomain of a cocartesian lifting of $f$ at the subobject. Dually, we think of the pre-image of a subobject of $Y$ as the domain of a cartesian lifting of $f$ at the suobject. This permits us to introduce and study a notion of topogenous orders on an appropriate faithful and amnestic functor which leads to the formal interior operators and includes the formal closure recently introdued as a particular case. While our motivation of studying topogenous structures on forms comes from topological functors, particularly the forgetful functor $F: \bf Top\longrightarrow \bf Set$, our approach permits to obtain topogenous orders which act on  quotients, subobjects,  as well as fibres of topological functors - see section 5 for some examples.
 
 {\indent In section 2, we recall a number of categorical concepts and results needed for the study of topogenous structures on a form. Section 3 studies the topogenous structures on forms. We show that formal closure operators form a reflective subcategory in the category of formal topogenous structures. Interior operators on forms are then introduced and, for forms in which every fibre is a complete lattice, these operators are shown to be in a one-to-one correspondence with a special class of formal topogenous orders. In section 4, we 
 study strict and final morphisms relative to the topogenous order introduced. Among other things, it is proved that a morphism is strict relative to a formal topogenous order if the order preserves the codomains of cocartesian liftings of the morphism. A final morphism relative to this order is the one for which the order reflects domains of its cartesian liftings. The note ends with section 5, which presents a number of examples that demonstrate our results.


\section{ Preliminaries} 
{\indent We consider a functor $F: \mathcal{A}\longrightarrow \mathcal{C}$. For an object $X$ in $\mathcal{C}$, a \textit{fibre} over $X$ is the subcategory $F^{-1}X$ of $\mathcal{A}$ consisting of those objects $A$ for which $FA = X$ and those morphisms $g$ satisfying $Fg = 1_{X}$ (see e.g \cite{MR2240597}). The functor $F$ is \textit{faithful} when for any $f, g: A\longrightarrow B$ in $\mathcal{A}$, $Ff = Fg \Rightarrow f = g$. $F$ is \textit{amnestic} if for any isomorphism $f: A\longrightarrow B$ in $\mathcal{A}$, $Ff = 1_{FA}\Rightarrow A = B$ and $f = 1_{A}$. For any $\mathcal{C}$-morphism $f : X\longrightarrow Y$, we define the relation $\leq_{f}: F^{-1}X\longrightarrow F^{-1}Y$ by : $A\leq_{f}  B$ if and only if there is a $\mathcal{A}$-morphism $A\longrightarrow B$ such that $F(A\longrightarrow B) = f$. In case $f = 1_{X}$, we write 
$A\leq_{X}  B$ (sometimes we omit the subscript and write $A\leq B$). If $A\leq B$ and $B\leq A$, then $A$ and $B$ are said to be \textit{fibre-isomorphic}. Of course, if $F$ is amnestic and $A$ and $B$ are fibre-isomorphic, then they are isomorphic (expressed as $A\cong B$). For any object $X\in \mathcal{C}$, $\leq_{X}$ is reflective and transitive. If  the functor $F$ is faithful, then $\leq_{X}$ is anti-symmetric if and only $F$ is amnestic.} 

{\indent Let $F: \mathcal{A}\longrightarrow \mathcal{C}$ be a faithful functor. According to \cite {MR2122803}, an $\mathcal{A}$-morphism $\alpha: A\longrightarrow B$ is an $F$-\textit{lifting} (or simply a lifting) of a $\mathcal{C}$-morphism  $f: FA\longrightarrow FB$ if $F\alpha = f$}. We may sometimes say that a $\mathcal{C}$-morphism  $f: FA\longrightarrow FB$ is an $\mathcal{A}$-morphism if it has a lifting $\alpha: A\longrightarrow B$. A \textit{cartesian lifting} of a $\mathcal{C}$-morphism  $f: X\longrightarrow Y$ at $B\in F^{-1}Y$ is a lifting of $f$ with condomain $B$ having the property that for any $C\in \mathcal{A}$, a $\mathcal{C}$-morphism  $g: FC\longrightarrow FA$ is an $\mathcal{A}$-morphism whenever  $f\circ g: FC\longrightarrow FB$ is an $\mathcal{A}$-morphism. Clearly, $\alpha: A\longrightarrow B$ is a cartesian lifting of $f: X\longrightarrow Y$ at $B\in F^{-1}Y$ if and only if $A\leq_{f}  B$ and for any $C\in \mathcal{A}$ and any $\mathcal{C}$-morphism  $g: FC\longrightarrow FA$, $C\leq_{f\circ g}  B\Leftrightarrow C\leq_{g}  A$. Dually, a \textit{cocartesian lifting} of $f: X\longrightarrow Y$ at $A\in F^{-1}X$ is a lifting of $f$ with domain $A$ having the property that for any $\mathcal{A}$-object $C$, a $\mathcal{C}$-morphism $g: FB\longrightarrow FC$ is an $\mathcal{A}$-morphism whenever $g\circ f : FA\longrightarrow FC$ is an $\mathcal{A}$-morphism. Clearly, $\alpha: A\longrightarrow B$ is a cocartesian lifting of $f: X\longrightarrow Y$ at $A\in F^{-1}X$ if and only if $A\leq_{f}  B$ and for any $C\in \mathcal{A}$ and any $\mathcal{C}$-morphism $g: FB\longrightarrow FC$,  $B\leq_{g}  C\Leftrightarrow A\leq_{g\circ f}  C$. The codomain of a cocartesian lifting of $f$ at $A$, when it exists, shall be denoted by $f.\;A$. Dually, the domain of a cartesian lifting of $f$ at $B$, when it exists, shall be denoted by $B.\;f$.  Thus, if $\alpha : B.\;f\longrightarrow B\;(\alpha :A\longrightarrow f.\;A)$ is a cartesian (co-cartesian) lifting of $f: X\longrightarrow Y$ at $B$ ($A$), then $C\leq_{g}  B.\;f\Leftrightarrow C\leq_{f\circ g}  B\;(f.\;A\leq_{g}  C\Leftrightarrow A\leq_{g\circ f} C)$ whenever $C\in \mathcal{A}$ and $g: FC\longrightarrow X\;(g: Y\longrightarrow FC)$ is a $\mathcal{C}$-morphism.} 
    
{\indent Following \cite {MR1313497}, a functor $F: \mathcal{A}\longrightarrow \mathcal{C}$ is a \textit{fibration}  if for every morphism $f: X\longrightarrow Y$ and $B\in F^{-1}Y$, there is a cartesian lift of $f$  at $B$. If  $F^{op}:\mathcal{A}^{op}\longrightarrow \mathcal{C}^{op}$, where the upper index "op" stands for "opposite",  is a fibration, then $F$ is said to be an \textit{opfibration}.
$F$ is $bifibration$ if it is both a fibration and an opfibration. In accordance with \cite{MR3226619}, by a \textit{form} over a category $\mathcal{C}$, we understand a faithfull and amnestic functor $F: \mathcal{A}\longrightarrow \mathcal{C}$. If $F$ is the codomain functor $ \mathcal{M}\longrightarrow \mathcal{C}$  where $\mathcal{M}$ is a class of monomorphisms in $\mathcal{C}$ seen as the full subcategory of the arrow category $\mathcal{C}^{2}$, then $F$ is called the  \textit{form of $\mathcal{M}$-subobjects}. The dual concept to the one of the form of $\mathcal{M}$-subobjects is the \textit{form of $\mathcal{M}$-quotients}.   According to \cite{MR3542097}, a form $F$ over a category $\mathcal{C}$ is said to be \textit{locally bounded} if each of its fibres has an upper bound and a lower bound. The upper bound (resp. the lower bound) of $F^{-1}X$ will be denoted by $1^{X}$ (resp. $0^{X}$). $F$ is said to be \textit{bounded} when it is locally
bounded and for any morphism $f: X\longrightarrow Y$ in $\mathcal{C}$, both $f.\;1^{X}$ and $0^{X}.\;f$ exist. The next Lemma that we recall from \cite{MR3542097} is a consequence of properties of (co)cartesian liftings of morphisms.

 \begin{Lemma}\label{L1}
  Let $F$ be a form over $\mathcal{C}$ and $f: X\longrightarrow Y$, $g: Y\longrightarrow Z$ $\mathcal{C}$-morphisms. 
  \begin{itemize}
\item [$(1)$] For any $A\in F^{-1}X$, if $f.\;A$ exists, then $g.\;(f.\;A)$ exists if and only if $(g\circ f).\;A$ exists, in which case $g.\;(f.\;A) = (g\circ f).\;A$. Dually, for any $B\in F^{-1}Z$, 
  if $B.\;g$ exists, then  $(B.\;g).\;f)$ is exists if and only if  $B.\;(g\circ f)$
exists, in which case $(B.\;g).\;f = B.\;(g\circ f)$.
\item [$(2)$] For any $A_{1}, A_{2}\in F^{-1}X$, if both $f.\;A_{1}$ and  $f.\;A_{2}$ exist, then $A_{1}\leq A_{2}$ implies that $f.\;A_{1}\leq f.\;A_{2}$. Dually, for any $B_{1}, B_{2}\in F^{-1}Y$, if both $B_{1}.\;f$ and $B_{2}.\;f$ exist, then $B_{1}\leq B_{2}$ implies that $B_{1}.\;f\leq B_{2}.\;f$ .
\end{itemize}
\end{Lemma}
 
\begin{defn}\label{D1}
\textit{A closure operator $C$} on a form $F: \mathcal{A}\longrightarrow \mathcal{C}$ or a \textit{formal closure operator} is a family of maps $\{C_{X}: F^{-1}X\longrightarrow \;F^{-1}X\;|\; X\in \mathcal{C}\}$ such that
\item [$(C1)$]  $A\leq C_{X}(A)$ for all $A\in \;F^{-1}X$ and $X\in \mathcal{C}$.
\item [$(C2)$]  For any morphism $f: X\longrightarrow Y$ in $\mathcal{C}$,  $A\leq_{f} B\Rightarrow C_{X}(A)\leq_{f} C_{Y}(B)$ for all $A\in \;F^{-1}X$
	      and $B\in \;F^{-1}Y$. 
\end{defn}
 We shall denote by Clo($F$) the conglomerate of all closure operators on $F$. Clo($F$) is ordered by $\preceq$ as follows: for any $C,C'\in \;$Clo($F$)  $C\preceq C'$ if and only if  $C_{X}(A)\leq C'_{X}(A)$ for all $A\in \;F^{-1}X$ and $X\in \mathcal{C}$. A closure operator $C$ on $F$ is idempotent if $C\circ C = C$, that is $C_{X}(C_{X}(A)) = C_{X}(A)$ for all $A\in \;F^{-1}X$ and $X\in \mathcal{C}$. If $F$ is a form such that for any $f: X\longrightarrow Y$ in $\mathcal{C}$, $A\in F^{-1}X$ and $B\in F^{-1}Y$, both $f.\;A$ and $B.\;f$ exist, then $(C2)$ in Definition \ref{D1} is equivalent to $(C3)$ $f.\;A\leq_{Y} B\Rightarrow f.\;C_{X}(A)\leq_{Y} C_{Y}(A)$ and to $(C4)$ $A\leq_{X} B.\;f\Rightarrow C_{X}(A)\leq_{X}  C_{Y}(B).\;f$. We can also obtain $(C2)$ as a conjunction of  $(C2')$ $A\leq_{X} B\Rightarrow C_{X}(A)\leq_{X} C_{Y}(B)$ and $(C2'')$ $f.\;C_{X}(A)\leq_{Y} C_{Y}(f.\;B)$ for suitable $A$, $B$.}

\section{Topogenous orders on forms}
{\indent For the rest of the paper, we work with a form $F$ over a category $\mathcal{C}$ such that for any $f: X\longrightarrow Y$ in $\mathcal{C}$, $A\in F^{-1}X$ and $B\in F^{-1}Y$, both $f.\;A$ and $B.\;f$ exist. We have the following useful Lemma.
 \begin{Lemma}\label{L2}
 Let $F$ be a form over a category $\mathcal{C}$, $f: X\longrightarrow Y$ be a $\mathcal{C}$-morphism, and let $A\in \;F^{-1}X$ and $B\in \;F^{-1}Y$. Then $A\leq_{X} (f.\;A).\;f$ and $f.\;(B.\;f)\leq_{Y} B$.  Moreover, if $F$ reflects sections and $f$ is a section, then $A\cong (f.\;A).\;f$. Dually, if $F$ reflects retractions and $f$ is a retraction, then $f.\;(B.\;f)\cong B$. If $F$ reflects isomorphisms and $f$ is an isomorphism with inverse $g$, then for any  $A\in \;F^{-1}X$, $f.\;A \cong A.\;g$. 
\end{Lemma}
\begin{proof}
 Let $ \beta :(f.\;A).\,f\longrightarrow f.\;A$ be a cartesian lifting of $f$ at $f.\;A$. Since $f = f\circ 1_{X}: FA\longrightarrow F(f.\;A)$, there is a morphism $h : A\longrightarrow (f.\;A).\;f$ with  $Fh =  1_{X}$. Hence $A\leq_{X} (f.\;A).\;f$. A dual argument shows that $f.\;(B.\;f)\leq_{Y} B$. Let $ \alpha :A\longrightarrow f.\;A$ be a cocartesian lifting of $f$ at $A$, $f$ be a section and let $F$ reflect sections. Since $F(\alpha) = F(\beta) = f$, we get that $F(\beta\circ h) = F(\beta)\circ F(h) = F(\alpha)\circ 1_{X} = F(\alpha)$. Since $F$ is faithful, we have that $\alpha = \beta\circ h$ and $\alpha$ is a section because $F$ reflects sections.  Therefore, $h$ is a section, i.e. there is a morphism $k : (f.\;A).\,f\longrightarrow A$ such that $k\circ h = 1_{A}$. Now, $Fk = Fk\circ 1_{X} = Fk\circ Fh = F(k\circ h) = F1_{A} = 1_{X}$. Consequently $ (f.\;A).\;f \leq_{X} A$. Thus, $A$ and $ (f.\;A).\;f$ are fibre-isomorphic and hence isomorphic. A dual raisoning proves that $f.\;(B.\;f)\cong B$. The last part of the proof follows from the fact that $f$ is an isomorphism if and only if $f$ is a retraction and a section. 
\end{proof}
}
It is not difficult to see from Lemma \ref{L2} that $f.\;-: F^{-1}X\longrightarrow F^{-1}Y$ which maps every $A\in \;F^{-1}X$ to $f.\;A$, and $- .\;f: F^{-1}Y\longrightarrow F^{-1}X$, which maps every $B\in \;F^{-1}Y$ to $A.\;f$, form a Galois connection, i.e. $f.\;A\leq_{Y} B\Leftrightarrow A\leq_{X} B.\;f$. Furthermore, the analysis of the proof of Lemma  \ref{L2} shows that the condition that $F$ reflects sections can be weakened to the condition that co-cartesian liftings of $f$ preserve sections. Dually, the condition that $F$ reflects retractions can be weakened to the condition that cartesian liftings of $f$ preserve retractions.

\indent{Let $F$ be the form of $\mathcal{M}$-subobjects and $\mathcal{C}$ a finetely $\mathcal{M}$-complete category endowed with a proper ($\mathcal{E}$, $\mathcal{M}$)-factorization structure for morphisms. Then for any $\mathcal{C}$-morphism $f: X\longrightarrow Y$, $F^{-1}X$ is the subobject lattice while $f.\;m$ (resp. $n.\;f$ ) is simply the image (resp. pre-image) of a subobject, for appropriate $m$ and $n$. In addition, if $\mathcal{C}$ has products of pairs so that sections in $\mathcal{C}$ belong to $\mathcal{M}$, then the conditions in Lemma \ref{L2} are satisfied.}
\begin{defn}\label{D2}
 $A\;topogenous\;order$ $\sqsubset$ on $F$ is a family  $\{\sqsubset_{X}\; |\;X\in \mathcal{C}\}$ of binary relations, each
	$\sqsubset_{X}$ on $F^{-1}X$, such that:
     \begin{itemize}
		\item [$(T1)$] $A\sqsubset_{X} B\Rightarrow A\leq_{X} B$ for every $A, B\in F^{-1}X$.
	    \item [$(T2)$]	$A'\leq_{X} A\sqsubset_{X} B\leq_{X} B'\Rightarrow A'\sqsubset_{X} B'$ for every $A, B\in F^{-1}X$.
		\item [$(T3)$] For every morphism $f: X\longrightarrow Y$ in $\mathcal{C}$, $f.\;A\sqsubset_{Y} B\Rightarrow A\sqsubset_{X} B.\;f$ for  $A\in F^{-1}X$ and $B\in F^{-1}Y$.
	\end{itemize}		
 \end{defn}	
 {\indent   It is quite clear that when $F$ is the form of $\mathcal{M}$-subobjects and $\mathcal{C}$ is finitely $\mathcal{M}$-complete, Definition \ref{D2} gives the categorical topogenous structures. 
 We shall denote by TORD($F$) the conglomerate of all topogenous structures on $F$. One orders TORD($F$) by $\subseteq $ as follows:  for any $\sqsubset, \sqsubset'\in $TORD($F$), $\sqsubset\subseteq \sqsubset'$ if and only if for all $X\in \mathcal{C}$ and $A, B\in \;F^{-1}X$, $A\sqsubset_{X} B\Rightarrow A\sqsubset'_{X} B$. A topogenous order $\sqsubset$ on $F$ is interpolative if $A\sqsubset_{X} B$ implies that there is $C\in F^{-1}X$ such that $A\sqsubset_{X} C\sqsubset_{X} B$. We denote by INTORD($F$) the class of all interpolative topogenous structures on $F$.
 Consider the following conditions for topogenous orders on $F$. }
 \begin{itemize}
		\item [$(TJ)$] For $\mathcal{A}_{I} = \{A_{i}\;|\;i\in I\}\subseteq F^{-1}X$, ($\forall i\in I,\;A_{i}\sqsubset_{X} B)\Rightarrow  \bigvee A_{i}\sqsubset_{X} B$ if $\bigvee A_{i}$ exists.
		\item [$(TM)$] For $\mathcal{B}_{I} = \{B_{i}\;|\;i\in I\}\subseteq_{X} F^{-1}X$, ($\forall i\in I,\;A\sqsubset_{X} B_{i})\Rightarrow  A\sqsubset_{X} \bigwedge B_{i}$ if $\bigwedge B_{i}$ exists.
 	\end{itemize}
 The conglomerate of all topogenous orders on $F$ satisfying condition $(TM)$ (resp.$(TJ)$)  will be denoted by MTORD($F$) (resp. JTORD($F$)). MTORD($F$) and JTORD($F$) are closed under arbitrary intersections in TORD($F$)  and are thus reflective subcategories. We next show that, for an appropriate form, the topogenous orders satisfying condition $(TM)$ are in one-to-one correspondence with the closure operators on $F$.
 \begin{Proposition}\label{P1}
Let $F$ be a form such that each of its fibres is a complete lattice.  For a topogenous order $\sqsubset$ on $F$ satisfying condition $(TM)$ and a closure operator $C$ on $F$, the assignments 
   $$ C\longmapsto \sqsubset^{C}\;\mbox{and}\;\sqsubset\longmapsto C^{\sqsubset}$$ where $A\sqsubset^{C}_{X} B\Leftrightarrow C_{X}(A)\leq_{X} B$ and $C^{\sqsubset} = \bigwedge \{B\in \;F^{-1}X\;|\;A\sqsubset_{X} B\}$, for any $X\in \mathcal{C}$and $A, B\in F^{-1}X$ define order isomorphisms inverse to each other between MTORD($F$) and Clo($F$). Moreover, $ C^{\sqsubset}$ is idempotent if and only if $\sqsubset$ is interpolative.
 \end{Proposition}
 \begin{proof}
 For  $\sqsubset^{C}$, $(T1)$ and  $(T2)$ are easily seen to be satisfied, and $(C1)$ is clear for $C^{\sqsubset}$.
 Let $f: X\longrightarrow Y$ be a $\mathcal{C}$-morphism and $\sqsubset \in $ TORD($F$). Then for any $A\in F^{-1}X$ and $B\in F^{-1}Y$ such taht $A\leq_{f} B$, $\{B'.\;f\;|\;B\sqsubset B'\}\subseteq \{A'\;|\;A\sqsubset A'\} $ by $(T3)$. This implies that $C^{\sqsubset}_{X}(A)\leq_{X}  C^{\sqsubset}_{Y}(B).\;f\Leftrightarrow C^{\sqsubset}_{X}(A)\leq_{f} C^{\sqsubset}_{Y}(B).$ For $(T3)$, $f.\;A\sqsubset_{Y} B\Leftrightarrow C^{\sqsubset}_{Y}(f.\;A)\leq_{Y} B\Rightarrow C^{\sqsubset}_{Y}(A)\leq_{f} B\Leftrightarrow C^{\sqsubset}_{Y}(A)\leq_{X} B.f\Leftrightarrow A\sqsubset_{X} B.\;f.$ Clearly, $C\longrightarrow \sqsubset^{C}\;\mbox{and}\;\sqsubset\longrightarrow C^{\sqsubset}$   preserve order and are inverse to each other. The fact that $ C^{\sqsubset}$ is idempotent if and only if $\sqsubset$ is interpolative is clear from the construction of $C^{\sqsubset} $.
 \end{proof}
Proposition \ref{P1} together with the well established relationship between categorical closure operators, interior operators and topogenous structures (see \cite{MR0286814, iragi2016quasi}) permit us to introduce formal interior operator and demonstrate that formal topogenous orders satisfying condition $(TJ)$ are indeed in a one-to-one correspondence with formal interior operators provided $F$ is a form such that each of its fibres is a complete lattice .  
\begin{defn}\label{D3}
An  \textit{operator} $I$ on a form $F: \mathcal{A}\longrightarrow \mathcal{C}$ or \textit{formal interior operator} is a family of maps $\{I_{X}: F^{-1}X\longrightarrow \;F^{-1}X\;|\; X\in \mathcal{C}\}$ such that
\item [$(I1)$]  $I_{X}(A)\leq A$ for all $A\in \;F^{-1}X$ and $X\in \mathcal{C}$.
\item [$(I1)$]   $A\leq B\Rightarrow I_{X}(A)\leq I_{X}(B)$ for all $A, B\in \;F^{-1}X$ and $X\in \mathcal{C}$.
\item [$(I3)$]  For any morphism $f: X\longrightarrow Y$ in $\mathcal{C}$,  $I_{Y}(B).\;f \leq I_{X}(B.\;f)$ for $B\in \;F^{-1}Y$. 
\end{defn}
We shall denote by Int($F$) the conglomerate of all interior operators on $F$. It is ordered by $I\leq I'$ if $I_{X}(A)\leq I'_{X}(A)$ for all $A\in \;F^{-1}X$ and $X\in \mathcal{C}$. An interior operator $I$ on $F$ is idempotent if $I\circ I = I$, that is $I_{X}(I_{X}(A)) = I_{X}(A)$ for all $A\in \;F^{-1}X$ and $X\in \mathcal{C}$.
A reasonning similar to the one in Proposition \ref{P1} results in the following proposition.

\begin{Proposition}\label{P2}
    Let $F$ be a form such that each of its fibres is a complete lattice. For a topogenous order $\sqsubset$ on $F$ satisfying condition $(TJ)$ and an interior operator $I$ on $F$, the assignments $$ I\longmapsto \sqsubset^{I}\;\mbox{and}\;\sqsubset\longmapsto I^{\sqsubset}$$ where $A\sqsubset^{I}_{X} B\Leftrightarrow A\leq_{X} I_{X}(A)$ and $I^{\sqsubset}(B) = \bigvee \{A\in F^{-1}X\;|\;A\sqsubset_{X} B\}$ for any $X\in \mathcal{C}$and $A, B\in F^{-1}X$ define order isomorphisms inverse to each other between ITORD($F$) and Int($F$). $ I^{\sqsubset}$ is idempotent if and only if $\sqsubset$ is interpolative.
  \end{Proposition} 
\section{Strict and final morphisms} 

\begin{Proposition}\label{P3}
Let $f: X\longrightarrow Y$ be a $\mathcal{C}$-morphism and $\sqsubset\in TORD($F$).$ Then $(T3)$ in Defintion \ref{D1} is equivalent to $A\sqsubset_{Y} B\Rightarrow A.\;f\sqsubset_{X} B.\;f$ for all $A,B\in F^{-1}Y$.
\end{Proposition}
Looking at axiom $(T3)$ in Definition \ref{D2} and Proposition \ref{P2}, one would ask the question
which morphisms satisfy the converse implication. We study in this section $\mathcal{C}$-morphisms that preserve  topogenous orders $\sqsubset$ on $F$ as well as those that reflect them.
\begin{defn}
Let $\sqsubset$ be a topogenous order on $F$. A morphism $f: X\longrightarrow Y$ in $\mathcal{C}$ is said to be 
$\sqsubset$-\textit{strict} if $A\sqsubset_{X} B.\;f\Rightarrow f.\;A\sqsubset_{Y} B$ for all $A\in F^{-1}X$ and $B\in F^{-1}Y$. Similarly, $f$ is $\sqsubset$-\textit{final} if $ B.\;f\sqsubset_{X} B'.\;f \Rightarrow B\sqsubset_{Y} B'$ for all   $B',B\in F^{-1}Y$.
\end{defn}
\begin{Theorem}
Let $\sqsubset$ be a topogenous order on $F$. A morphism $f: X\longrightarrow Y$ is $\sqsubset$-stict if and only if $\sqsubset$ preserves the codomains of cocartesian liftings of $f$,  i.e. $A\sqsubset_{X} B\Rightarrow f.\;A\sqsubset_{X} f.\;B$ for  $A, B\in F^{-1}X$.
\end{Theorem}
\begin{proof}
Assume that $f$ is $\sqsubset$-strict and $A\sqsubset_{X} B$. Since cartesian and cocartesian liftings of $f$ exist, $(f.\;B).\;f$ exists and $A\sqsubset_{X} B\leq (f.\;B).\;f\Rightarrow A\sqsubset_{X} (f.\;B).\;f\Rightarrow f.\;A\sqsubset_{X} f.\;B$. Conversely, if $\sqsubset$ preserves codomains of cocartesian liftings of $f$ and $ A\sqsubset_{X} B.\;f $, then $f.\;(B.\;f)$ exists and $ A\sqsubset_{X} B.\;f\Rightarrow f.\;A\sqsubset_{X} f.\;(B.\;f)\leq B\Rightarrow f.\;A\sqsubset_{X} B$.
\end{proof}
 Let us recall from \cite{MR3542097} the notion of a thick morphism which will help us to characterize  $\sqsubset$-final morphisms. In a bounded form $F$ over $\mathcal{C}$,  $f: X\longrightarrow Y\in \mathcal{C}$ is a \textit{thick morphism} if $f.\;1^{X} = 1^{Y}.$

 \begin{Theorem}
 Let $F$ be a locally bounded form and $\sqsubset \in $TORD($F$). Then, every $\sqsubset$-final morphism $f: X\longrightarrow Y$ is thick provided $1^{Y}\sqsubset 1^{Y}.$ The condition $1^{Y}\sqsubset 1^{Y}$ can be dropped if $\sqsubset \in $MTORD($F$).  If $F$ reflects retractions, then every retraction which is $\sqsubset$-final  is $\sqsubset$-strict while if $F$ reflects sections, every section that is $\sqsubset$-strict is $\sqsubset$-final. 
 \end{Theorem}
 \begin{proof}
  Since $1^{Y}\sqsubset 1^{Y}$ and $f$ is $\sqsubset$-final, $1^{Y}\sqsubset 1^{Y}\Rightarrow 1^{Y}.\;f\sqsubset 1^{Y}.\;f\Rightarrow 1^{Y}.\;f\sqsubset 1^{X}\leq (f.\;1^{X}).\;f\Rightarrow 1^{Y}.\;f\sqsubset 1^{X}\leq (f.\;1^{X}).\;f\Rightarrow 1^{Y}\sqsubset f.\;1^{X}\Rightarrow 1^{Y}\leq f.\;1^{X}\Rightarrow f.\;1^{X} = 1^{Y}$. It is clear that if $\sqsubset \in $MTORD($F$), then $1^{Y}\sqsubset 1^{Y}$ can be dropped because $1^{Y}\sqsubset 1^{Y}\Leftrightarrow C^{\sqsubset}(1^{Y}) = 1^{Y}.$ Assume $F$ reflects retractions and $f$ is a retraction that is $\sqsubset$-final.  Then, by Lemma \ref{L2}, $A\sqsubset_{X} B\Leftrightarrow (f.\;A).\;f\sqsubset (f.\;B).\;f\Leftrightarrow f.\;A\sqsubset f.\;B$. Lastly, let $F$ reflect sections and $f$ be a section which is $\sqsubset$-strict. Then, by Lemma \ref{L2}, $ B.\;f\sqsubset_{X} B'.\;f \Rightarrow f.\;(A.\;f)\sqsubset f.\;(B.\;f)\Rightarrow A\sqsubset B.$
  
  \end{proof}
 \begin{Proposition}
 Let $\sqsubset$ be a topogenous order on $F$. 
  If $F$ reflects isomorphisms, then the class of $\sqsubset$-strict (resp. $\sqsubset$-final) morphisms contains all isomorphisms of $\mathcal{C}$. The class of $\sqsubset$-strict (resp. $\sqsubset$-final) morphisms is closed under composition.  If $F$ reflects retractions, $g\circ f$ is a $\sqsubset$-strict (resp. $\sqsubset$-final) morphism, and $f$ is a retraction, then $g$ is $\sqsubset$-strict (resp. $\sqsubset$-final) morphism. Dually, if $F$ reflects sections, $g\circ f$ is a $\sqsubset$-strict (resp. $\sqsubset$-final) morphism, and $g$ is a retraction,  $f$ is $\sqsubset$-strict (resp. $\sqsubset$-final).
 \end{Proposition}
 \begin{proof}
Assume that $f$ is an isomorphism with inverse $g$. Since $F$ reflects isomorphisms, by Lemmas \ref{L2} and \ref{L1}, we have that $A\sqsubset B.\;f\Rightarrow f.\;A = A.\;g \sqsubset (B.\;f).\;g = B.\;(g\circ f) = B.\;1_{X} = B$ for any  $A\in \;F^{-1}X$ and  $B\in \;F^{-1}Y$.  If $f: X\longrightarrow Y$ and $g: Y\longrightarrow Z$ are $\sqsubset$-strict, then by Lemma \ref{L1},  $A\sqsubset B.\;(f\circ g) = (B.\;f).\;g \Leftrightarrow f.\;A\sqsubset B.\;g \Leftrightarrow (g\circ f).\;A = g.\;(f.\;A)\sqsubset B.$
Let $f$ be a retraction and $g\circ f$ $\sqsubset$-strict. Since $F$ reflects retractions, By Lemmas \ref{L2} and \ref{L1}, $A\sqsubset B.\;g\Rightarrow A.\;f\sqsubset (B.\;g).\;f = B.\;(g\circ f)\Rightarrow (g\circ f).(A.\;f)\sqsubset B\Leftrightarrow g.\;[f.\;(A.\;f)]\sqsubset B\Leftrightarrow g.\;A\sqsubset B.$
Let $f$ be a section and $g\circ f$ $\sqsubset$-strict. Since $F$ reflects sections, by Lemmas \ref{L2} and \ref{L1}, $A\sqsubset B.\;f = [(g.\;A).\;g].\;f = (g.\;B).\;(g\circ f)\Rightarrow g.\;(f.\;A)\sqsubset g.\;B\Rightarrow f.\;A\sqsubset (g.\;B).\;g\Rightarrow f.\;A\sqsubset B$. A similar resonning can be applied for the case of $\sqsubset$-final.
  
 \end{proof}
{\indent Taking into consideration Propositions \ref{P2} and \ref{P1} and the fact that for any morphism $f: X\longrightarrow Y$ in $\mathcal{C}$ cocartesian (resp. cartesian) liftings of $f$ at any $A\in F^{-1}X$ (resp. $B\in F^{-1}Y$) exist, we obtain the following result.}
\begin{Proposition}
Let $\sqsubset$ be a topogenous order on $F$, $f: X\longrightarrow Y$ a $\mathcal{C}$-morphism and assume that $\sqsubset\in $JTORD($F$). Then $f$ is $\sqsubset$-stict if $I^{\sqsubset}$  preserves domains of cartesian liftings of $f$, i.e. $ I^{\sqsubset}(B).\;f = I^{\sqsubset}_{X}(B.\;f)$ for any $B\in F^{-1}Y$. Similarly, if $\sqsubset\in $MTORD($F$), then $f$ is $\sqsubset$-strict if and only if $C^{\sqsubset}$ preserves codomains of cocartesian liftings of $f$, i.e. $f.\;C^{\sqsubset}_{X}(A) = C^{\sqsubset}_{X}(f.\;A)$ for any $A\in F^{-1}X$.
\end{Proposition}
\indent{ Let $\mathcal{C}$ be finitely cocomplete (so that pushouts of $\mathcal{E}$-morphisms along arbirary $\mathcal{C}$-morphisms exist and are in $\mathcal{E}$)  with a proper $(\mathcal{E}, \mathcal{M})$-factorization structure for morphisms. Let $F$ be the domain functor $dom: \mathcal{E}\longrightarrow \mathcal{C}$. For any  $X\in \mathcal{C}$, $F^{-1}X$ is the preordered class of $\mathcal{E}$-quotients of $X$. For any $\mathcal{C}$-morphism $f: X\longrightarrow Y$ and $e\in F^{-1}X$, $f.\;e$} is the pushout of $e$ along $f$ while $d .\;f$ is the $\mathcal{E}$-part of the $(\mathcal{E}, \mathcal{M})$-factorization of $d\circ f$ for any $d\in F^{-1}Y$}.
\begin{defn}
A topogenous order $\sqsubset$ on $F$ is said to be \textit{cohereditary} if for any retraction  $f: X\longrightarrow Y$ and any $A, B\in F^{-1}Y$, $ B.\;f\sqsubset_{X} A.\;f \Rightarrow B\sqsubset_{Y} A$. Equivalently, $\sqsubset$ is said to be \textit{cohereditary} if every retraction is $\sqsubset$-final. 
\end{defn}
\indent{When $\sqsubset\in $JTORD($F$), $\sqsubset$ is cohereditary if and only if $C^{\sqsubset}_{Y}(B) = f.\;C^{\sqsubset}_{X}(B.\;f)$ for any retraction in $\mathcal{C}$. It was observed in \cite{MR3641248} that there is an order reversing isomorphism between the poset of full $\mathcal{E}$-reflective replete subcategories of $\mathcal{C}$ and the poset of cohereditary idempotent closure operators on the form of $\mathcal{E}$-quotient. This result together with Proposition \ref{P1} permit to affirm that \textit{there is an order reversing isomorphism between the poset of full $\mathcal{E}$-reflective replete subcategories of $\mathcal{C}$ and the conglomerate of all cohereditary and interpolative topogenous orders on the form of $\mathcal{E}$-quotients in $\mathcal{C}$ satisfying condition (TM)}. }


\section{Some examples}

 $\bf I.\;$ Let $F$ be the forgetful functor from $\bf Top$ to $\bf Set$.  For every  $X\in \bf Set$, $F^{-1}X$ is the complete lattice of all topologies on $X$. If $f : X \longrightarrow Y$  is a function in $\bf Set$, $f.\;\mathcal{T}_{X}$ is the final topology on $Y$ induced by $f$ for all $\mathcal{T}_{X}\in F^{-1}X,$ and $\mathcal{T}_{Y}.\;f$ is the initial topology on $X$ induced by $f$ for all $\mathcal{T}_{Y}\in F^{-1}Y$. $F$ is a bounded form in which $\mathcal{T}_{Y}.\;f$ and $f.\;\mathcal{T}_{X}$ exist for any $\mathcal{T}_{X}\in F^{-1}X$ and $\mathcal{T}_{Y}\in F^{-1}Y$.  Let $\mathcal{T}_{X}, \mathcal{T'}_{X}\in F^{-1}X.$

 {\indent $(a)$ Putting $\mathcal{T}_{X}\sqsubset_{X} \mathcal{T'}_{X}\Leftrightarrow \theta(\mathcal{T}_{X})\leq_{X} \mathcal{T'}_{X}$,  where $\theta(\mathcal{T}_{X})$ is the $\theta$-topology generated by $\mathcal{T}_{X}$, we get a topogenous order on $F$. $(T1)$ and $(T2)$ follow from the fact that $\theta(\mathcal{T}_{X})\leq_{X} \mathcal{T'}_{X}\Leftrightarrow \mathcal{T'}_{X}\subseteq \theta(\mathcal{T}_{X})$. For $(T3)$, let $f: X\longrightarrow Y$ and $\mathcal{T}_{Y},\mathcal{T'}_{Y} \in F^{-1}Y$. Assume that $\mathcal{T'}_{Y}\subseteq \theta(\mathcal{T}_{Y})$ and $ A\in \mathcal{T'}_{Y}.\;f$. Then there is $O'\in \mathcal{T'}_{Y}$ such that $A = f^{-1}(O')$ and $O'\in \theta(\mathcal{T}_{Y})$ by the assumption. Let $y\in A$. Then $f(y)\in O$ and so there is $U\in \mathcal{U}_{f(y)}$ with $U$ closed and $U\subseteq O$. This implies that $f^{-1}(U)\subseteq f^{-1}(O) = A$ and $f^{-1}(U)\in \mathcal{U}_{Y}$. Consequently $A\in \theta(\mathcal{T}_{Y}.\;f)$, i.e. $\mathcal{T'}_{Y}.\;f\subseteq\theta(\mathcal{T}_{Y}.\;f) $. It is easy to see that $\sqsubset \in MTORD(F)$. }
 
 \begin{Proposition}
A function $f: X\longrightarrow Y$ is  $\sqsubset$-strict if it is surjective and for any $\mathcal{T}_{X}\in F^{-1}X$ and $\mathcal{T}_{Y}\in F^{-1}Y$, $f: (X,\mathcal{T}_{X})\longrightarrow (Y,\mathcal{T}_{Y})$ is clopen.
Every surjective function $f: X\longrightarrow Y$ is $\sqsubset$- final.
 \end{Proposition}
\begin{proof}
 Assume $\mathcal{T}_{Y}.\;f\subseteq \theta(\mathcal{T}_{X})$ and $A\in \mathcal{T}_{Y}$. If $y\in A$, then by surjectivity of $f$, there is $x\in f^{-1}(A)$ such that $f(x) = y$ and  $f^{-1}(A)\in \mathcal{T}_{Y}.\;f$. By assumption, there are $O\in \mathcal{T}_{X}$ and $U$ closed in $\mathcal{T}_{X}$ such that $x\in O\subseteq U\subseteq f^{-1}(A)$. Since $f$ is clopen, $f(O)\in \mathcal{T}_{Y}$ and $f(U)$ is closed in $\mathcal{T}_{Y}$. 
 We have $f(x) = y\in f(O)\subseteq f(U)\in A$. Thus $A\in \theta(f.\;\mathcal{T}_{X})$, that is  $\mathcal{T}_{Y}\subseteq \theta(f.\;\mathcal{T}_{X})$. Consequently $f$ is $\sqsubset$-strict. Assume that $\mathcal{T'}_{Y}.\;f\subseteq \theta(\mathcal{T}_{Y}.\;f)$ and $A\in \mathcal{T'}_{Y}.$
 Let $x\in A$. Since $f$ is surjective, there is $x\in f^{-1}(A)$ such that $y = f(x)$ and $f^{-1}(A)\in \mathcal{T'}_{Y}.\;f$. By the assymption, $f^{-1}(A)\in \theta(\mathcal{T}_{Y}.\;f)$ and so there are open $O\in \mathcal{T}_{Y}.\;f$ and closed $U$ in $\mathcal{T}_{Y}.\;f$ such that $x\in O\subseteq U\subseteq f^{-1}(A)$. Now,
 $O = f^{-1}(O')$ and $U = f^{-1}(U')$ with $O'\in \mathcal{T}_{Y}$ and $U'$ closed in $\mathcal{T}_{Y}$. We get that $f(x)\in O'\subseteq U'\subseteq A$. Thus $A\in \theta(\mathcal{T}_{Y})$ and $\mathcal{T'}_{Y}\subseteq \theta(\mathcal{T}_{Y}).$
 \end{proof}
 \indent{$(b)$ Putting $\mathcal{T}_{X}\sqsubset_{X} \mathcal{T'}_{X}\Leftrightarrow \mathcal{T}_{X}\leq_{X} b(\mathcal{T'}_{X})$, where $b(\mathcal{T'}_{X})$ is the $b$-topology generated by $\mathcal{T'}_{X}$, is a topogenous order on $F$. $(T1)$ and $(T2)$ follow from the fact that $\mathcal{T}_{X})\leq_{X} b(\mathcal{T'}_{X})\Leftrightarrow b(\mathcal{T'}_{X})\subseteq \mathcal{T}_{X}$. For $(T3)$, let $f: X\longrightarrow Y\in \bf Set$ and $\mathcal{T}_{X}\in F^{-1}X$, $\mathcal{T}_{Y}\in F^{-1}Y$. Assume that 
 $b(\mathcal{T}_{Y})\subseteq f.\;\mathcal{T}_{X}$ and $A\in \mathcal{T}_{Y}.\;f.$ Then, there is  $O\in \mathcal{T}_{Y}$ such that $A = f^{-1}(O)$. Since $\mathcal{T}_{Y}\subseteq b(\mathcal{T}_{Y})$, $O\in b(\mathcal{T}_{Y})$ and by the assumption, $O\in f.\;\mathcal{T}_{X}$ which implies that $A\in \mathcal{T}_{X}$ and $A\in b(\mathcal{T}_{X})$. Thus $\mathcal{T}_{Y}.\;f\subseteq b(\mathcal{T}_{X})$. Clearly, $\sqsubset \in JTORD(F)$.  
 \begin{Proposition}
 Every function $f: X\longrightarrow Y$ is $\sqsubset$ strict (resp. $\sqsubset$-final). 
 \end{Proposition}
 \begin{proof}
  Assume that $b(\mathcal{T'}_{X})\subseteq \mathcal{T}_{X}$ and $A\in b(f.\mathcal{T'}_{X})$. If $x\in f^{-1}(A)$, then $f(x)\in A$ and there are $O\in f.\;\mathcal{T'}_{X}$ and $F\in f.\;\mathcal{T}_{X}$ such that $f(x)\in O\cap F\subseteq A.$ Now, $f^{-1}(O)\in \mathcal{T'}_{X}$ and $f^{-1}(F)$ is closed in $\mathcal{T'}_{X}$. We have that $x\in f^{-1}(O\cap F) =  f^{-1}(O)\cap  f^{-1}(F)) \subseteq f^{-1}(A)$. Since $f^{-1}(O)\in \mathcal{T'}_{X}$  and $f^{-1}(F)$ is closed in $\mathcal{T'}_{X}$, $f^{-1}(A)\in b(\mathcal{T'}_{X})$ which implies that $f^{-1}(A)\in \mathcal{T}_{X}$. Thus $A\in b(f.\mathcal{T}_{X})$.
  Assume $b(\mathcal{T'}_{Y}.\;f)\subseteq \mathcal{T}_{Y}.\;f$ for any $\mathcal{T}_{Y}, \mathcal{T'}_{Y}\in F^{-1}Y$  and $A\in b(\mathcal{T'}_{Y})$.  Let $x\in f^{-1}(A)$. Then $f(x)\in A$ and there are $O\in \mathcal{T'}_{Y}$ and $F$ closed in $\mathcal{T'}_{Y}$ such that  $f(x)\in O\cap F\subseteq A$. This implies that  $x\in f^{-1}(O)\cap f^{-1}(F)\subseteq f^{-1}(A)$. Thus $f^{-1}(A)\in b(\mathcal{T'}_{Y}.\;f)\subseteq \mathcal{T}_{Y}.\;f$.
 \end{proof}
   }
{\indent   $\bf II.$ Another example, somehow similar to the previous one, is the concrete functor $F$ from $\bf Qunif$ to $\bf Top$ (where $\bf Top$ denotes the category of topological spaces and continuous maps and  $\bf Qunif$) denotes the category of quasi-uniform spaces and quasi-uniformly continuous maps. Then, for any $X\in \bf Top$, $F^{-1}X$ is the complete lattice of all quasi-uniform structures on $X$ compatible with $\mathcal{T}$, known as \textit{functorial quasi-uniform structures}. If $f : X \longrightarrow Y$ is a continuous map and $\mathcal{U}_{Y}\in F^{-1}Y$, then $\mathcal{U}_{Y}.\;f$ is initial quasi-uniformity induced by $f$ and $f.\;\mathcal{U}_{X}$ is the largest quasi-uniformity on $Y$ for which $f$ is quasi-uniformly continuous. For any $\mathcal{U}_{X}, \mathcal{V}_{X}\in F^{-1}X$, putting $\mathcal{U}_{X}\sqsubset_{X} \mathcal{V}_{X}\Leftrightarrow \mathcal{U}_{X}\leq_{X} \mathcal{V}^{\star}_{X}$, where $\mathcal{V}^{\star}_{X}$ is the coarsest uniformity containing $\mathcal{V}_{X}$, we get a topogenous order on $F$.} 

 \bigskip
 
{\indent  $\bf III.$ Let $F$ the form of subgroups (i.e injective group homomorphisms) with $\bf Grp$ the category of groups and group homomorphisms. Then, for any $X\in \bf Grp$, $F^{-1}X$ is the complete lattice of all subgroups of $X$. It is also clear that this form is bounded and  $A.\;f$ and $f.\;B$ exist for any $A\in F^{-1}X$ and $B\in F^{-1}Y.$ Now, define $\sqsubset$ on $F$ by $A\sqsubset B\Leftrightarrow A\leq N\leq B$ with $N$ a normal subgroup of $X$. Then $\sqsubset$ is a topogenous order on $F$. A group homomorphism $f : X \longrightarrow Y$ is $\sqsubset$-strict if and only if it preserves normal subgroups while $f$ is $\sqsubset$-final if and only if it is surjective.}

\bibliographystyle{abbrv}
\bibliography{references}
\addcontentsline{toc}{chapter}{References}
\end{document}